\documentclass[a4paper,reqno,12pt]{amsart}
\usepackage{amssymb}

\usepackage{amsfonts}
\usepackage{amsmath}
\usepackage{cite}
\usepackage{amssymb}
\usepackage{mathrsfs}
\usepackage{hyperref}
\numberwithin{equation}{section}
\newtheorem{theorem}{Theorem}[section]
\newtheorem{prop}[theorem]{Proposition}
\newtheorem{lm}[theorem]{Lemma}

\newcommand{\R}{\mathbb{R}}

\newcommand{\N}{\mathbb{N}}

\def\ds{\displaystyle}

\def\dive{\mathrm{div}\ }
\newenvironment{preuve}{{\noindent {\bf Proof. }}}{\hfill {\rule{2.5mm}{2.5mm}} }

\def\1{\textbf{1\!\!1}}
\newenvironment{preuve1.3}{{\noindent {\bf Proof of Theorem \ref{thm1.3}. }}}{\hfill {\rule{2.5mm}{2.5mm}} }

\newenvironment{preuve1.4}{{\noindent {\bf Proof of Theorem \ref{thm1.4}. }}}{\hfill {\rule{2.5mm}{2.5mm}} }

\author[M.~Amara]{Mustapha Amara}
\author[J.~Benameur]{Jamel Benameur}
\address{University of Gabes, Faculty of Science of Gab\`es, Department of Mathematics, Research Laboratory Mathematics and Applications LR17ES11; Tunisia}
\email{\sl jamelbenameur@gmail.com}
\email{\sl Mostafa.Amara@fsg.u-gabes.tn}

\title[Global solution of anisotropic Quasi-Geostrophic Equations in Sobolev Space]
{Global solution of anisotropic Quasi-Geostrophic Equations}

\begin{document}
\begin{abstract}
	In \cite{YZ}, the author proved the global existence of the two-dimensional anisotropic quasi-geostrophic equations with condition on the parameters $\alpha,$ $\beta$  in the Sobolev spaces $H^s( \R^2)$; $s\geq 2$. In this paper, we show that this equations has a global solution in the spaces $H^s(\R^2)$, where $\max\{2-2\alpha,2-2\beta\}< s<2$, with additional condition over $\alpha$ and $\beta$. The proof is based on the Gevrey-class regularity of the solution in neighborhood of zero.
\end{abstract}



\maketitle
\tableofcontents


\section{\bf Introduction}
	We interest to the study the two-dimensional surface quasi-geostrophic  equation $(SQG)$ with fractional horizontal dissipation	and fractional vertical thermal diffusion, which can be written as
	$$
		(AQG)	\left\{\begin{array}{l}
			\partial_t\theta+ u_\theta.\nabla\theta +(\mu|\partial_1|^{2\alpha}+\nu |\partial_2|^{2\beta})\theta=0,x=(x_1,x_2)\in\R^2, \ t>0,\\
			u_\theta=\mathcal{R}^\perp \theta,\\
			\theta(x,0)=\theta^0(x),\ x\in \R^2.
		\end{array}\right.
	$$
where $u_\theta=(u_\theta^1,u_\theta^2)$ is the velocity is determined by the Riesz transforms of the potential temperature $\theta$ via the formula
\begin{align*}
	&u_\theta^1=-\partial_{2}(-\Delta)^{-\frac{1}{2}}\theta,
	&u_\theta^2= \partial_{1} (-\Delta)^{-\frac{1}{2}}\theta,
\end{align*}
and $\alpha,\beta\in (0,1)$ and $\mu,\nu> 0$ are real
		constants.
\vskip0.5cm
\par The $(AQG)$ equation is an important model of geophysical fluid dynamics, which describes the evolution of the surface temperature field in the rotating stratified fluid. The first mathematical studies of the $(AQG)$ equation was carried out in 1994s by Constantin,
Majda and Tabak. For more details and mathematical and physical explanations of this model we can consult \cite{CP}, \cite{CP1}, \cite{DC} and \cite{JP}.\\

In addition to the physical interpretation of the $(AQG)$ equation, it can also be used as a simplified model of the 3D Navier-Stokes equation.\\

Easy to see that when $\alpha=\beta$ and $\mu =\nu$, the dissipation term $\mu|\partial_{1}|^{2\alpha}\theta +\nu|\partial_{2}|^{2\beta}\theta $ can be considered as the standard fractional Laplacian $(-\Delta)^{\alpha}\theta$, $\left(\left(|\xi_1|^{2\alpha}+|\xi_2|^{2\alpha}\right)\sim |\xi|^{2\alpha}\right)$, therefore $(AQG)$ becomes the classical dissipation $(QG)$ equation, with its form as follows
$$(QG)
	\left\{\begin{array}{l}
		\partial_t\theta+\mu(-\Delta)^{\alpha}\theta+ u_\theta.\nabla\theta =0,x=(x_1,x_2)\in\R^2, \ t>0\\
		\theta(x,0)=\theta^0(x).
	\end{array}\right.
$$
The local well-posedness of $(QG)$ with $H^{2-2\alpha}(\R^2)$ data is established by  \cite{JN} .\\

The global regularity results of the system $(AQG)$ for the case $\alpha,\beta\in (0,1)$ satisfy
\begin{equation}\label{P}
	\beta>\begin{cases}
		\frac{1}{2\alpha+1},&0<\alpha\leq \frac{1}{2}\\
		\\
		\frac{1-\alpha}{2\alpha},& \frac{1}{2}<\alpha<1.	
	\end{cases}
\end{equation}
were considered by Zhuan Ye, \cite{YZ}, in $H^s(\R^2)$ for $s\geq 2$. The condition (\ref{P}) can be written as follows:
$$Y=\{(\alpha,\beta)\in(0,1)^2:\;(\alpha,\beta)\;satisfies\;(\ref{P})\}=Y_1\cup Y_2\cup Y_3,$$
where
$$\begin{array}{lcl}
Y_1&=&\{(\alpha,\beta)\in(0,1)^2:\;\alpha,\beta\in(1/2,1)\;and\;\beta>\frac{1-\alpha}{2\alpha}\}\\
&=&(1/2,1)^2\\
Y_2&=&\{(\alpha,\beta)\in(0,1)^2:\;\alpha\in(1/2,1)\;and\;\beta\in(\frac{1-\alpha}{2\alpha},1/2]\}\\
Y_3&=&\{(\alpha,\beta)\in(0,1)^2:\;\alpha\in(0,1/2]\;and\;\beta\in(\frac{1}{2\alpha+1},1)\}.
\end{array}$$

 In this paper we study the case $(\alpha,\beta)\in Y_1=(1/2,1)^2$. We start by proving the local well-posedness of the equations $(AQG) $ in $H^s(\R^2)$ for $s>\max\{2-2\alpha,2-2\beta\}$, where we have applied the Fixed Point Theorem to the integral form of our system. The second result is to prove the Gevrey-class regularity of this solution near to zero:
	$$\left(t\mapsto e^{t(|\partial_{1}|^\alpha+|\partial_{2}|^{\beta})}\theta\right)\in L^\infty\left([0,\varepsilon],H^s(\R^2)\right),$$
which implies that $\theta\in C((0,\varepsilon),H^2(\R^2)$. By combining this result with \cite{YZ}, we obtain the global solution.

In this paper, all constants will be denoted by $C$ that is a generic constant depending only on the quantities specified in the context.
\section{\bf Main theorems}
To simplify the calculus and some steps of proofs of our results, we assume that $\mu=\nu = 1$, then we get the following main theorems:

We move now for the main result of the paper,
\begin{theorem} \label{thm1.4}
	Let $\alpha,\beta\in (1/2,1)$ and $\theta^0\in H^s(\R^2)$, for $s\in (\max\{2-2\alpha,2-2\beta\},2).$ Then, the system $(AQG)$ admits a unique global solution
	$$\theta\in C(\R^+,H^s(\R^2));\;|\partial_1|^\alpha\theta,|\partial_2|^\beta\theta\in L^2_{loc}(\R^+,H^s(\R^2)).$$
\end{theorem}

For the proof of the global existence of our equation, we need Zhuan Ye's theorem in \cite{YZ} as follows:
\begin{theorem}\label{thmYE}\cite{YZ}
	Let $\theta^0\in H^s(\R^2)$ with $s\geq 2$ and $\alpha, \beta>0$. Then there exists a positive $T(\|\theta^0\|_{H^s})>0$ such that for the system $(AQG)$
	admits a unique solution $$\theta\in C([0,T],H^s(\R^2)),\quad|\partial_{1}|^\alpha\theta,|\partial_{2}|^\beta\theta\in L^2([0,T],H^s(\R^2)).$$
	Moreover, if $\alpha,\beta\in (0,1)$ satisfy \eqref{P}, then, the system $(AQG)$ admits a unique global solution $\theta$ such that for any $T>0$ we have
	$$
		\theta\in C([0,T],H^s(\R^2)),\ |\partial_{1}|^{\alpha}\theta,|\partial_{2}|^{\beta}\theta\in L^2([0,T],H^s(\R^2)).
	$$
\end{theorem}

 	\section{\bf Notation and Preliminary Results}
 In this short section, we collect some notations and definitions that will be used later, and we	give some technical lemmas.
 \begin{itemize}
 	\item[$\bullet$] The Fourier transformation in $\R^2$		
 	$$
 		\mathcal{F}(f)(\xi)=\widehat{f}(\xi)=\int_{\R^2}e^{-ix.\xi}f(x)dx,\quad \xi\in \R^2.
 	$$
 	The inverse Fourier formula is
 	$$
 		\mathcal{F}^{-1}(f)(x)=(2\pi)^{-2}\widehat{f}(\xi)=\int_{\R^2}e^{i\xi.x}f(\xi)d\xi,\quad x\in \R^2.
 	$$
 	\item[$\bullet$] The fractional operators
 	$\ds|\partial_1|=\sqrt{-\partial_{x_1}^2}$ and $\ds|\partial_2|=\sqrt{-\partial_{x_2}^2}$ are defined through the Fourier transform, namely
 	$$
 		|\partial_1|^{2\alpha}f=\mathcal{F}^{-1}\left(\xi\mapsto|\xi_1|^{2\alpha}\widehat{f}(\xi)\right),\quad |\partial_2|^{2\beta}f=\mathcal{F}^{-1}\left(\xi\mapsto|\xi_2|^{2\beta}\widehat{f}(\xi)\right).
 	$$
 	\item[$\bullet$] The convolution product of a suitable pair of function $f$ and $g$ on $\R^2$ 	is given by
 	$$
 		f\ast g(x)=\int_{\R^2} f(x-y)g(y)dy
 	$$
 	\item[$\bullet$] If $f=(f_1,f_2)$ and $g=(g_1,g_2)$ are two vector fields, we set
 	$$f\otimes g:= (g_1f,g_2f),$$
 	and
 	$$\dive (f\otimes g):= (\dive(g_1f),\dive(g_2f)).$$
 	\item[$\bullet$] Let $(X, \|.\|_{X})$, be a Banach space, $p\in [0,+\infty]$ and $T > 0$. We define $$L_T^p(X):=L^p([0,T],X)$$ the space of all
 	measurable functions $t\in[0,T]\mapsto f(t)\in X $ such that
 	$$\left(t\mapsto \|f(t)\|_{X}\right)\in L^p([0,T]).$$
 	\item[$\bullet$] For $s\in \R$, the space:
 	$$H^{s}(\R^2):=\left\{f\in \mathcal{S}'(\R^2); (1+|\xi|^2)^{s/2}\widehat{f}\in L^2(\R^2)\right\}$$ denotes the usual inhomogeneous Sobolev space on $\R^2$, with the norm
 	$$\|f\|_{H^s}=\left(\int_{\R^2}(1+|\xi|^2)^s|\widehat{f}(\xi)|^2d\xi\right)^{\frac{1}{2}}.$$
 \item[$\bullet$] For $s\in \R$, the space:
 	$$\dot{H}^{s}(\R^2):=\left\{f\in \mathcal{S}'(\R^2);\widehat{f}\in L^1(\R^2)\mbox{ and }|\xi|^s\widehat{f}\in L^2(\R^2)\right\}$$ denotes the usual homogeneous Sobolev space on $\R^2$, with the norm
 	$$\|f\|_{\dot{H}^s}=\left(\int_{\R^2}|\xi|^{2s}|\widehat{f}(\xi)|^2d\xi\right)^{\frac{1}{2}}.$$
 	\item[$\bullet$] For $s_1,s_2\in\R$ and $t\in[0,1]$; the interpolation inequalities, respectively, in the homogeneous and non-homogeneous Sobolev spaces
 	\begin{align}
 		&\|f\|_{H^{ts_1+(1-t)s_2}}\leq \|f\|_{H^{s_1}}^t\|f\|_{H^{s_2}}^{1-t},\\
 		&\|f\|_{\dot{H}^{ts_1+(1-t)s_2}}\leq \|f\|_{\dot{H}^{s_1}}^t\|f\|_{\dot{H}^{s_2}}^{1-t}.
 	\end{align}
 \item For $s\in\R$, $a>0$ and $\sigma>1$, the Gevrey-Sobolev space is defined as follows
 $$H^s_{a,\sigma}(\R^2)=\{f\in {\mathcal S}'(\R^2);\;e^{a|\xi|^{1/\sigma}}(1+|\xi|^2)^{s/2}\widehat{f}\in L^2(\R^2)\}.$$
 \end{itemize}

 We recall a fundamental lemma concerning the Sobolev injection and some product laws in homogeneous Sobolev
 spaces.
 \begin{lm}\cite{BH}\label{lemma1.2}
 	Let $p \in [2, +\infty)$ and $\sigma\in [0,1)$ such that
 	$$\frac{1}{p}+\frac{\sigma}{2}=\frac{1}{2}.$$
 	Then, there is a constant $C > 0$ such	 that
 	$$\|f\|_{L^p(\R^d)}\leq C\||\nabla|^{\sigma}f\|_{L^2(\R^d)}.$$
 \end{lm}
 \begin{lm}\cite{BH}\label{lmprod}
 	Let $s_1$, $s_2$ be two real numbers such that $s_1<1$ and $s_1+s_2>0$. Then, there exists a positive constant $C=C(s_1,s_2)$ such that for all $f,g\in \dot{H}^{s_1}(\R^2)\bigcap \dot{H}^{s_2}(\R^2)$;
 	\begin{equation}
 		\|fg\|_{\dot{H}^{s_1+s_2-1}}\leq C(s_1,s_2) \left(\|f\|_{\dot{H}^{s_1}}\|g\|_{\dot{H}^{s_2}}+\|f\|_{\dot{H}^{s_2}}\|g\|_{\dot{H}^{s_1}}\right).
 	\end{equation}
 	Moreover, in addition $s_2<1$, there exists a positive constant $C'=C'(s_1,s_2)$ such that, for all $f\in \dot{H}^{s_1}(\R^2)$ and $g\in \dot{H}^{s_2}(\R^2)$;
 	\begin{equation}
 		\|fg\|_{\dot{H}^{s_1+s_2-1}}\leq C'(s_1,s_2)\|f\|_{\dot{H}^{s_1}}\|g\|_{\dot{H}^{s_2}}
 	\end{equation}
 \end{lm}
 \begin{lm}[\cite{JN} Calderòn-Zygmund Theorem]\label{lemma1.1}
 	For any $p\in (1,+\infty)$, there is a constant $C(p)>0$ such that
 	\begin{equation}
 		\|\mathcal{R}^\perp\theta\|_{L^p}\leq C(p) \|\theta\|_{L^p}.
 	\end{equation}
 \end{lm}
\par We recall the following important commutator and product estimates: The following lemma due to Sobolev's interpolation theorem
	\begin{lm}\label{lm3.4}
	Let $s\in \R$ and $\alpha,\beta\in (0,1)$, such that $\alpha\leq \beta$ then for any $f\in \mathcal{S}$
	\begin{equation*}
		\||\nabla|^\alpha f\|_{\dot{H}^s}\leq \| f\|_{\dot{H}^s} + \||\partial_{1}|^\alpha f\|_{\dot{H}^s}+ \||\partial_{2}|^\beta f\|_{\dot{H}^s}.
	\end{equation*}
\end{lm}
\begin{preuve}
	We have
	\begin{align*}
		\||\nabla|^\alpha f\|_{\dot{H}^s}^2&=\int_{\R^2}(|\xi_1|^2+|\xi_2|^2)^{\alpha}\left|\mathcal{F}(|\nabla|^{s}f)(\xi)\right|^2d\xi\\
		&\leq \int_{\R^2}|\xi_1|^{2\alpha}\left|\mathcal{F}(|\nabla|^{s}f)(\xi)\right|^2d\xi+\int_{\R^2}|\xi_2|^{2\alpha}\left|\mathcal{F}(|\nabla|^{s}f)(\xi)\right|^2d\xi=\||\partial_{1}|^\alpha f\|_{\dot{H}^s}^2+\||\partial_{2}|^\alpha f\|_{\dot{H}^s}^2.
	\end{align*}
	But $0<\alpha\leq \beta$ then their exist $z\in (0,1]$ such that $\alpha=z\times\beta+(1-z)\times0$ and by interpolation inequality we get
	\begin{align*}
		\||\partial_{2}|^\alpha f\|_{\dot{H}^s}&\leq \| f\|_{\dot{H}^s}^{1-z}\||\partial_{2}|^\beta f\|_{\dot{H}^s}^{z}\\
		&\leq (1-z)\| f\|_{\dot{H}^s}+z\||\partial_{2}|^\beta f\|_{\dot{H}^s} \leq \| f\|_{\dot{H}^s}+\||\partial_{2}|^\beta f\|_{\dot{H}^s}.
	\end{align*}
\end{preuve}
\begin{lm}\cite{BJ}\label{lm3.8}
 	Let $d\in \N$ and $r\in (0,1].$ Then
 	\begin{equation}
 		|\xi|^r\leq |\xi-\eta|^r+|\eta|^r,\quad \forall \xi,\eta\in \R^d.
 	\end{equation}
 \end{lm}
 \begin{lm}\label{lmexp}
 	Let $\alpha>0$ and $x,r>0$, there is a constant $C(\alpha)>0$ such that
 	\begin{equation}
 		x^\alpha e^{-rx}\leq \frac{C(\alpha)}{r^\alpha},\;\forall x>0.
 	\end{equation}
 \end{lm}
 \begin{preuve}
 	We consider the following function
 	$$\begin{array}{cccc}
 		f:&[0,+\infty)&\rightarrow&\R\\
 		&x&\mapsto&f(x)=x\exp({-x})
 	\end{array}.$$
 	We have, for any $x>0$, $f(x)=x\exp({-x})<1$. So
 	\begin{align*}
 		x^\alpha e^{-rx} &=\frac{\alpha^\alpha}{r^\alpha}\left(\frac{r}{\alpha}x \exp\left({-\frac{r}{\alpha}x}\right)\right)^\alpha\\
 		&\leq \frac{\alpha^\alpha}{r^\alpha}.
 	\end{align*}
 	Which implies the result.
 \end{preuve}

  	\section{\bf Proof of Theorem \ref{thm1.3}}
  	The case of $s\geq 2$ is treated in \cite{YZ}, then we want to prove the result for $s\in (2-2\alpha,2)$. This proof is done in four steps: in the first, we prove the uniqueness and local existence $\theta\in C([0,T_0],H^s(\R^2))$. The second step is devoted to prove the Gevrey-class regularity near to 0; $e^{t(|\partial_1|^\alpha+|\partial_2|^\beta)}\theta\in C([0,T_1],H^s(\R^2))$,\;$0<T_1\leq T_0$. In the third step, we prove that $\theta\in C((0,T_1),H^2(\R^2))$. In the last step and by combining the previous steps and Theorem \ref{thmYE}, we prove the global existence in $H^s$.
\subsection{Step 1: Uniqueness and local existence} For $T>0$ we define the closed subset of $ C([0,T],H^{s}(\R^2))$
  	$$\mathcal{B}(T)=\left\{\theta\in C([0,T],H^{s}(\R^2)):\ \|\theta\|_{L^\infty_T(H^{s})}\leq 2\|\theta^0\|_{H^s}\right\}.$$
  	Now, Consider the following application
  	$$\begin{array}{cccl}
  		\psi:&\mathcal{B}(T)&\rightarrow&C([0,T],H^{s}(\R^2))\\
  		&\theta&\mapsto&\psi(\theta)=L_0-B(\theta,\theta)
  	\end{array}$$
  	where
  	$$\begin{cases}
  		L_0(t)=e^{-tA(D)}\theta^0=\mathcal{F}^{-1}\left(\xi\mapsto e^{-tA(\xi)}\widehat{\theta^0}(\xi)\right),\\
  		\ds B(\theta_1,\theta_2)(t)=\int_{0}^{t}e^{-(t-\tau)A(D)}\dive(\theta_1  u_{\theta_2})d\tau.
  	\end{cases}$$
  	where $A(\xi)=|\xi_1|^{2\alpha}+|\xi_2|^{2\beta}$. We will apply the fixed point theorem with a good choice of  $T$: \\
  	We are looking for a condition on $ T$ such that $\psi(\mathcal{B}(T))\subset \mathcal{B}(T)$. For that we have
  	$$	\|L_0(t)\|_{H^s}=\left(\int_{\R^2}|\xi|^{2s}e^{-2tA(\xi)}|\widehat{\theta^0}(\xi)|^2d\xi\right)^{\frac{1}{2}}\leq \left(\int_{\R^2}|\xi|^{2s}|\widehat{\theta^0}(\xi)|^2d\xi\right)^{\frac{1}{2}} =\|\theta^0\|_{H^{s}}.$$
  	For the estimation of $B(\theta_1,\theta_2),$ we have two cases; $s\in (2-2\alpha,1)$ and $s\in[1,2)$.\\
  	We start with the first case, when $s\in (2-2\alpha,1)$ and  we have  for any $\theta_1,\theta_2 \in C([0,T],H^s(\R^2))$:
  	\begin{align*}
  		\|B(\theta_1,\theta_2)(t)\|_{H^s}&\leq \|B(\theta_1,\theta_2)(t)\|_{L^2}+\|B(\theta_1,\theta_2)(t)\|_{\dot{H}^s}\\
  		&\leq \int_{0}^{t} \|e^{-(t-\tau)A(D)}\dive(\theta_1u_{\theta_2})\|_{L^2}d\tau+\int_{0}^{t} \|e^{-(t-\tau)A(D)}\dive(\theta_1u_{\theta_2})\|_{\dot{H}^s}d\tau.
  	\end{align*}
  	We must estimate each term alone, we see the term with $\dot{H}^s$ norm and we can write
  	\begin{align*}
  		 \|e^{-(t-\tau)A(D)}\dive(\theta_1u_{\theta_2})\|^2_{\dot{H}^s}&\leq\int_{\R^2}|\xi|^{2(s+1)}e^{-2(t-\tau)A(\xi)}|\widehat{\theta_1u_{\theta_2}}(\xi)|^2d\xi\\
  		&\leq\int_{\R^2}|\xi|^{4-2s}e^{-2(t-\tau)A(\xi)}|\xi|^{2(2s-1)}|\widehat{\theta_1u_{\theta_2}}(\xi)|^2d\xi.
  	\end{align*}
  	Using Lemma \ref{lmexp}, we get
  	\begin{align*}
  		|\xi|^{4-2s}e^{-2(t-\tau)A(\xi)}&\leq |\xi_1|^{4-2s}e^{-(t-\tau)|\xi_1|^{2\alpha}}+|\xi_2|^{4-2s}e^{-(t-\tau)|\xi_2|^{2\beta}}\\
  		&\leq  C \left(\frac{1}{(t-\tau)^{\frac{2-s}{\alpha}}}+ \frac{1}{(t-\tau)^{\frac{2-s}{\beta}}}\right).
  	\end{align*}
  	Then, due of Lemma \ref{lmprod} with $s_1=s_2=s$ we obtain
  	\begin{align*}
  		\|e^{-(t-\tau)A(D)}\dive(\theta_1u_{\theta_2})\|_{\dot{H}^s}
  		&\leq C \left(\frac{1}{(t-\tau)^{\frac{2-s}{2\alpha}}}+\frac{1}{(t-\tau)^{\frac{2-s}{2\beta}}}\right)\left\|\theta_1u_{\theta_2}\right\|_{\dot{H}^{2s-1}}\\
  		&\leq C\left(\frac{1}{(t-\tau)^{\frac{2-s}{2\alpha}}}+\frac{1}{(t-\tau)^{\frac{2-s}{2\beta}}}\right)\left\|\theta_1\right\|_{H^s}\left\|\theta_2\right\|_{H^s}.
  	\end{align*}
  	Moreover, for the second term with $L^2$ norm, we have
  	\begin{align*}
  		\|e^{-(t-\tau)A(D)}\dive(\theta_1u_{\theta_2})\|_{L^2}&=\left(\int_{\R^2}|\xi|^{2}e^{-(t-\tau) A(\xi)}|\widehat{\theta_1u_{\theta_2}}(\xi)|^2d\xi\right)^{\frac{1}{2}}\\
  		&=\left(\int_{\R^2}|\xi|^{4-2s}e^{-(t-\tau) A(\xi)}|\xi|^{2(s-1)}|\widehat{\theta_1u_{\theta_2}}(\xi)|^2d\xi\right)^{\frac{1}{2}}\\
  		&\leq C \left(\frac{1}{(t-\tau)^{\frac{2-s}{2\alpha}}}+\frac{1}{(t-\tau)^{\frac{2-s}{2\beta}}}\right)\left\|\theta_1u_{\theta_2}\right\|_{\dot{H}^{s-1}}\\
  		&\leq C \left(\frac{1}{(t-\tau)^{\frac{2-s}{2\alpha}}}+\frac{1}{(t-\tau)^{\frac{2-s}{2\beta}}}\right)\left\|\theta_1\right\|_{H^s}\left\|\theta_2\right\|_{H^s},
  	\end{align*}
  	where we use Lemma \ref{lmprod} for $s_1=s$ and $s_2=0$. Therefore
  	\begin{align}
  		\nonumber	\|B(\theta_1,\theta_2)(t)\|_{H^s}&\leq  C \int_{0}^{t} \left(\frac{1}{(t-\tau)^{\frac{2-s}{2\alpha}}}+\frac{1}{(t-\tau)^{\frac{2-s}{2\beta}}}\right)\left\|\theta_1\right\|_{H^s}\left\|\theta_2\right\|_{H^s}d\tau\\
  		\label{B(u,v)}	&\leq  C_1(\alpha,\beta,s)  \left(T^{\frac{s-2+2\alpha}{2\alpha}}+T^{\frac{s-2+2\beta}{2\beta}}\right) \left\|\theta_1\right\|_{L^\infty_T(H^s)}\left\|\theta_2\right\|_{L^\infty_T(H^s)}.
  	\end{align}
  	As a result, we have for any $\theta\in \mathcal{B}(T)$:
  	\begin{align*}
  		\|\psi(\theta)\|_{H^s}&\leq \|L_0(t)\|_{H^{s}}+\|B(\theta,\theta)(t)\|_{H^s}\\
  		&\leq \|\theta^0\|_{H^s}+	C_1(\alpha,\beta,s)  \left(T^{\frac{s-2+2\alpha}{2\alpha}}+T^{\frac{s-2+2\beta}{2\beta}}\right) \left\|\theta\right\|_{L^\infty_T(H^s)}^2\\
  		&\leq \|\theta^0\|_{H^s}+ 4\	C_1(\alpha,\beta,s)  \left(T^{\frac{s-2+2\alpha}{2\alpha}}+T^{\frac{s-2+2\beta}{2\beta}}\right) \left\|\theta^0\right\|^2_{H^s}.
  	\end{align*}
  	By choosing $T>0$ such that
  	$$  T^{\frac{s-2+2\alpha}{2\alpha}}+T^{\frac{s-2+2\beta}{2\beta}}\leq  \frac{1}{8C_1(\alpha,\beta,s)\|\theta^0\|_{H^{s}}}.$$
  	We get $\psi (\mathcal{B}(T))\subset \mathcal{B}(T)$.\\
  	Looking now for a second condition on $ T$ such that $\psi$ is contracting on $\mathcal{B}(T)$, so, for any $\theta_2$, $\theta_2\in \mathcal{B}(T)$, and we have
  	\begin{align*}
  		\|\psi(\theta_1)-\psi(\theta_2)\|_{H^s}&\leq  \|B(\theta_1,\theta_1)-B(\theta_2,\theta_2)\|_{H^s}\\
  		&\leq \|B(\theta_1-\theta_2,\theta_1)+B(\theta_2,\theta_1-\theta_2)\|_{H^s}\\
  		&\leq \|B(\theta_1-\theta_2,\theta_1)\|_{H^s}+\|B(\theta_2,\theta_1-\theta_2)\|_{H^s}\\
  		&\leq C_1(\alpha,\beta,s) \left(T^{\frac{s-2+2\alpha}{2\alpha}}+T^{\frac{s-2+2\beta}{2\beta}}\right)\left(\|\theta_1\|_{L^\infty_T(H^s)}+\|\theta_2\|_{L^\infty_T(H^s)}\right) \|\theta_1-\theta_2\|_{L^\infty_T(H^s)}\\
  		&\leq 4C_1(\alpha,\beta,s) \|\theta^0\|_{H^{s}} \left(T^{\frac{s-2+2\alpha}{2\alpha}}+T^{\frac{s-2+2\beta}{2\beta}}\right) \|\theta_1-\theta_2\|_{L^\infty_T(H^s)},
  	\end{align*}
  	where we use the estimation \eqref{B(u,v)}. Then, for $T>0$ satisfy
  	\begin{equation}\label{T0}
  		T^{\frac{s-2+2\alpha}{2\alpha}}+T^{\frac{s-2+2\beta}{2\beta}}\leq \frac{1}{8C_1(\alpha,\beta,s)\|\theta^0\|_{H^{s}}},
  	\end{equation}
  	we ensure also that $\psi$ is contracting on $\mathcal{B}(T)$.\\
  	
  	Now we pass to the second case: if $1\leq s<2$. Let $\theta_1,\theta_2 \in C([0,T],H^s(\R^2))$, then, we have
  	\begin{align*}
  		\|B(\theta_1,\theta_2)(t)\|_{H^s}&\leq \int_{0}^{t} \left(\int_{\R^2}|\xi|^2e^{-(t-\tau)A(\xi)}|\widehat{\theta_1u_{\theta_2}}(\xi)|^2d\xi\right)^{\frac{1}{2}}d\tau\\
  		&\hskip0.5cm+\int_{0}^{t} \left(\int_{\R^2}|\xi|^{2(s+1)}e^{-(t-\tau)A(\xi)}|\widehat{\theta_1u_{\theta_2}}(\xi)|^2d\xi\right)^{\frac{1}{2}}d\tau.
  	\end{align*}
  	However, by Lemma \ref{lmexp} and Lemma \ref{lmprod} with $s_1=s_2=\frac{1}{2}$, we obtain
  	\begin{align}\label{3.3}
  		\nonumber\int_{0}^t	\left(\int_{\R^2}|\xi|^2e^{-(t-\tau)A(\xi)}|\widehat{\theta_1u_{\theta_2}}(\xi)|^2d\xi\right)^{\frac{1}{2}}	d\tau&\leq C
  		\int_{0}^t	\left(\frac{1}{(t-\tau)^\frac{1}{2\alpha}}+\frac{1}{(t-\tau)^\frac{1}{2\beta}}\right)\|\theta_1u_{\theta_2}\|_{L^{2}}	d\tau\\
  	\nonumber	&\leq C
  		\int_{0}^t	 \left(\frac{1}{(t-\tau)^\frac{1}{2\alpha}}+\frac{1}{(t-\tau)^\frac{1}{2\beta}}\right)\|\theta_1\|_{\dot{H}^{\frac{1}{2}}}\|u_{\theta_2}\|_{\dot{H}^{\frac{1}{2}}}	 d\tau\\
  		&\leq C(\alpha,\beta,s) \left(T^{\frac{2\alpha-1}{2\alpha}}+T^{\frac{2\beta-1}{2\beta}}\right)\|\theta_1\|_{L^\infty_T(H^s)} \| \theta_2\|_{L^\infty_T(H^s)}.
  	\end{align}
  	For the other term we have
  	\begin{align*}
  		 \left(\int_{\R^2}|\xi|^{2(s+1)}e^{-(t-\tau)A(\xi)}|\widehat{\theta_1u_{\theta_2}}(\xi)|^2d\xi\right)^\frac{1}{2}&=\left(\int_{\R^2}|\xi|^{2(1+\varepsilon)}e^{-(t-\tau)A(\xi)}|\xi|^{2(s-\varepsilon)}|\widehat{\theta_1u_{\theta_2}}(\xi)|^2d\xi\right)^\frac{1}{2}\\
  		&\leq C\left(\frac{1}{(t-\tau )^{\frac{1+\varepsilon}{2\alpha}}}+\frac{1}{(t-\tau )^{\frac{1+\varepsilon}{2\beta}}}\right) \|\theta_1u_{\theta_2}\|_{\dot{H}^{s-\varepsilon}},
  	\end{align*}
  	where $0<\varepsilon<2\alpha-1$ (we can choose $\varepsilon=\varepsilon(\alpha)=\frac{2\alpha-1}{2}$). Therefore, by Lemma \ref{lmprod} for $s_1=s$ and $s_2=1-\varepsilon$, we get
  	\begin{align}
  		\nonumber\int_{0}^t	\left(\int_{\R^2}|\xi|^{2(s+1)}e^{-(t-\tau)A(\xi)}|\widehat{\theta_1u_{\theta_2}}(\xi)|^2d\xi\right)^\frac{1}{2}d\tau
  		&\leq C\int_{0}^t\left(\frac{1}{(t-\tau )^{\frac{1+\varepsilon}{2\alpha}}}+\frac{1}{(t-\tau )^{\frac{1+\varepsilon}{2\beta}}}\right) \|\theta_1u_{\theta_2}\|_{\dot{H}^{s-\varepsilon}}d\tau\\
  		\nonumber	&\hskip-3.5cm\leq C\int_{0}^t\left(\frac{1}{(t-\tau )^{\frac{1+\varepsilon}{2\alpha}}}+\frac{1}{(t-\tau )^{\frac{1+\varepsilon}{2\beta}}}\right) \left(\|\theta_1\|_{\dot{H}^s}\| \theta_2\|_{\dot{H}^{1-\varepsilon}}+\|\theta_1\|_{\dot{H}^{1-\varepsilon}}\| \theta_2\|_{\dot{H}^{s}}\right)d\tau\\
  		\label{3.4}	&\hskip-3.5cm\leq C(\alpha,\beta,s)\left(T^{\frac{2\alpha-1-\varepsilon}{2\alpha}}+T^{\frac{2\beta-1-\varepsilon}{2\beta}}\right) \|\theta_1\|_{L^\infty_T(H^s)}\| \theta_2\|_{L^\infty_T(H^s)}
  	\end{align}
  	Then, by \eqref{3.3} and \eqref{3.4}, we have
  	\begin{align*}
  		\|B(\theta_1,\theta_2)(t)\|_{H^s}
  		&	\leq C_2(\alpha,\beta,s)\left(T^{\frac{2\alpha-1}{2\alpha}}+T^{\frac{2\beta-1}{2\beta}}+T^{\frac{2\alpha-1}{4\alpha}}+T^{\frac{4\beta-2\alpha-1}{4\beta}}\right)\| \theta_1\|_{L^\infty_T(H^{s})}\|\theta_2\|_{L^\infty_T(H^{s})}.
  	\end{align*}
  	Thus, if $\theta\in \mathcal{B}(T)$
  	\begin{align*}
  		\|\psi(\theta)\|_{H^s}&\leq \|L_0(t)\|_{H^{s}}+\|B(\theta,\theta)(t)\|_{H^s}\\
  		&\leq \|\theta^0\|_{H^s}+ 4\	C_2(\alpha,\beta,s)  \left(T^{\frac{2\alpha-1}{2\alpha}}+T^{\frac{2\beta-1}{2\beta}}+T^{\frac{2\alpha-1}{4\alpha}}+T^{\frac{4\beta-2\alpha-1}{4\beta}}\right) \left\|\theta^0\right\|^2_{H^s}.
  	\end{align*}
  	By choosing $T>0$ such that
  	\begin{equation}
  		T^{\frac{2\alpha-1}{2\alpha}}+T^{\frac{2\beta-1}{2\beta}}+T^{\frac{2\alpha-1}{4\alpha}}+T^{\frac{4\beta-2\alpha-1}{4\beta}}\leq  \frac{1}{8C_2(\alpha,\beta,s)\|\theta^0\|_{H^{s}}},
  	\end{equation}
  	we get $\psi (\mathcal{B}(T))\subset \mathcal{B}(T)$.\\
  	Now, we want to prove that $\psi$ is contracting by the best choose of $T$, so for $\theta_1$, $\theta_2\in \mathcal{B}(T)$, we have
  	\begin{align*}
  		\|\psi(\theta_1)-\psi(\theta_2)\|_{H^s}
  		&\leq \|B(\theta_1-\theta_2,\theta_1)\|_{H^s}+\|B(\theta_2,\theta_1-\theta_2)\|_{H^s}\\
  		&\hskip-2cm\leq C_2(\alpha,\beta,s)  \left(T^{\frac{2\alpha-1}{2\alpha}}+T^{\frac{2\beta-1}{2\beta}}+T^{\frac{2\alpha-1}{4\alpha}}+T^{\frac{4\beta-2\alpha-1}{4\beta}}\right)\left( \|\theta_1\|_{L^\infty_T(H^s)}+\|\theta_2\|_{L^\infty_T(H^s)}\right) \|\theta_1-\theta_2\|_{L^\infty_T(H^s)}\\
  		&\hskip-2cm\leq 4C_2(\alpha,\beta,s) \|\theta^0\|_{H^{s}} \left(T^{\frac{2\alpha-1}{2\alpha}}+T^{\frac{2\beta-1}{2\beta}}+T^{\frac{2\alpha-1}{4\alpha}}+T^{\frac{4\beta-2\alpha-1}{4\beta}}\right) \|\theta_1-\theta_2\|_{L^\infty_T(H^s)}.
  	\end{align*}
  	Then, for $T>0$ satisfy
  	\begin{equation}\label{T1}
  		T^{\frac{2\alpha-1}{2\alpha}}+T^{\frac{2\beta-1}{2\beta}}+T^{\frac{2\alpha-1}{4\alpha}}+T^{\frac{4\beta-2\alpha-1}{4\beta}}\leq \frac{1}{8C_2(\alpha,\beta,s)\|\theta^0\|_{H^{s}}},
  	\end{equation}
  	we ensure also that $\psi$ is contracting on $\mathcal{B}(T)$.\\
  	
  	Therefore, if $T_0>0$ satisfy
  	\begin{align}\label{T1}
  		&	T_0^{\frac{s-2+2\alpha}{2\alpha}}+T_0^{\frac{s-2+2\beta}{2\beta}}\leq \frac{1}{8C_1(\alpha,\beta,s)\|\theta^0\|_{H^{s}}},\\
  		&	T_0^{\frac{2\alpha-1}{2\alpha}}+T_0^{\frac{2\beta-1}{2\beta}}+T_0^{\frac{2\alpha-1}{4\alpha}}+T_0^{\frac{4\beta-2\alpha-1}{4\beta}}\leq \frac{1}{8C_2(\alpha,\beta,s)\|\theta^0\|_{H^{s}}},
  	\end{align}
  	then, we get $\psi(\mathcal{B}(T_0))\subset\mathcal{B}(T_0)$ and $\psi $ is contraction on $\mathcal{B}(T_0)$, for any $s\in(2-2\alpha,2)$.\\
  	Fixed Point Theorem gives the existence and uniqueness of solution of $(AQG)$ in $C([0,T_0],H^s(\R^2))$.\\
  	
\subsection{Step 2: Gevrey-class regularity of local solution} It is a question of proving the following proposition. \begin{prop}\label{thm1.3} Let $\alpha,\beta\in (1/2,1)$
	and $\theta^0\in H^{s}(\R^2)$, $\max\{2-2\alpha,2-2\beta\}<s<2$.
	Then, there exists
	a positive time $0<T_1\leq T_0$, such that the unique solution
	$\theta\in C([0,T_0];H^{s}(\R^2))$ given by the first step, satisfies
	\begin{equation}\label{1.2}
		\left(t\mapsto e^{t\left(|\partial_{1}|^{\alpha}+|\partial_{2}|^{\beta}\right)}\theta (t)\right)\in L^\infty([0,T_1];H^{s}(\R^2)).
	\end{equation}
\end{prop}
{\bf Proof.} The idea of proof is to apply the Fixed Point Theorem to a subpart of
${\mathcal B(T_0)}\subset C([0,T_0],H^s(\R^2))$ by using the same theory and the same steps of the preceding proof for the function $\Theta=e^{\frac{t}{2}B(D)}\theta(t)$.\\
  	Let $T>0$ and consider
  	$$\mathcal{P}(T)=\left\{\theta\in \mathcal{B}(T):\ \|e^{\frac{.}{2}B(D)}\theta\|_{L^\infty_T(H^s)}\leq 2\|\theta^0\|_{H^s}\right\},$$
  	where $e^{\frac{t}{2}B(D)}\theta=\mathcal{F}^{-1}\left(\xi\mapsto e^{\frac{t}{2}B(\xi)}\widehat{\theta}(\xi) \right)$ and $B(\xi)=2\left(|\xi_1|^{\alpha}+|\xi_2|^\beta\right)$.\\
  	
  	So we need to prove that $\psi (\mathcal{P}(T))\subset \mathcal{P}(T)$ and $\psi$  is contracting on $\mathcal{P}(T)$  for a best choose of $T$.\\
  	
  	Let $\theta\in \mathcal{P}(T)$, then for any $t\in [0,T]$, we have
  	$$|\widehat{\psi(\theta)}(t)|\leq e^{-tA(\xi)}|\widehat{\theta^0}(\xi)|+\int_{0}^{t}|\xi|e^{-(t-\tau)A(\xi)}|\widehat{\theta u_\theta}|d\tau.$$
  	We multiple it with  $(1+|\xi|^{2})^{\frac{s}{2}}e^{\frac{t}{2} B(\xi)}$, we get
  	\begin{align*} (1+|\xi|^{2})^{\frac{s}{2}}e^{\frac{t}{2} B(\xi)}|\widehat{\psi(\theta)}(t)|&\leq (1+|\xi|^{2})^{\frac{s}{2}} e^{-\frac{t}{2} (A(\xi)-B(\xi))}|\widehat{\theta^0}(\xi)|\\
  		&\hskip0.6cm+C\int_{0}^{t}|\xi|e^{-\frac{1}{2}(t-\tau) A(\xi)}e^{-\frac{1}{2}(t-\tau) (A(\xi)-B(\xi))}e^{\frac{\tau}{2}B(\xi)}|\widehat{\theta u_\theta}|d\tau\\
  		&\hskip0.6cm+C\int_{0}^{t}|\xi|^{s+1}e^{-\frac{1}{2}(t-\tau) A(\xi)}e^{-\frac{1}{2}(t-\tau) (A(\xi)-B(\xi))}e^{\frac{\tau}{2}B(\xi)}|\widehat{\theta u_\theta}|d\tau
  	\end{align*}
  	By Lemma \ref{lm3.8}, we get for any $\xi,\eta\in \R^2$
  	\begin{align*}
  		B(\xi)&=2(|\xi_1|^\alpha+|\xi_2|^{\beta})\\
  		&=2(|\xi_1-\eta_1+\eta_1|^\alpha+|\xi_2-\eta_2+\eta_2|^{\beta})\\
  		&\leq 2(|\xi_1-\eta_1|^\alpha+|\eta_1|^\alpha+|\xi_2-\eta_2|^\beta+|\eta_2|^{\beta})= B(\xi-\eta)+B(\eta),
  	\end{align*}
  which implies
   \begin{align*}
  	e^{\frac{\tau}{2}B(\xi)}|\widehat{\theta u_\theta}|&\leq\int_{\R^2}e^{\frac{\tau}{2} B(\xi-\eta)}e^{\frac{\tau}{2} B(\eta)} |\widehat{\theta}(\xi-\eta)||\widehat{\theta}(\eta)|d\xi=\left|\widehat{\Theta.\Theta}\right|,
  \end{align*}
  where $\Theta=\mathcal{F}^{-1}\left(\xi\mapsto e^{\frac{t}{2} B(\xi)} |\widehat{\theta}(\xi)|\right)$.\\
  Moreover, by using the elementary inequality, we get
  	\begin{align*}
  	A(\xi)-B(\xi)&=(|\xi_1|^{\alpha}+1)^2+(|\xi_1|^{\alpha}+1)^2-2\geq -2.
  \end{align*}
  	We deduce that
  	\begin{align*} \|e^{\frac{t}{2} B(D)}\psi(\theta)\|_{H^s}&\leq  e^{t}\|\theta^0\|_{H^s}+C\int_{0}^{t}e^{t-\tau}\left(\int_{\R^2}|\xi|^2e^{-(t-\tau) A(\xi)}|\widehat{\Theta.\Theta}|^2d\xi\right)^{\frac{1}{2}}d\tau\\
  		&\hskip0.6cm+C\int_{0}^{t}e^{t-\tau}\left(\int_{\R^2}|\xi|^{2(s+1)}e^{-(t-\tau) A(\xi)}|\widehat{\Theta.\Theta}|^2d\xi\right)^{\frac{1}{2}}d\tau\\
  		&\leq e^{T}\|\theta^0\|_{H^s}+Ce^{T}\int_{0}^{t}\left(\int_{\R^2}|\xi|^2e^{-(t-\tau) A(\xi)}|\widehat{\Theta.\Theta}|^2d\xi\right)^{\frac{1}{2}}d\tau\\
  		&\hskip0.6cm+Ce^{T}\int_{0}^{t}\left(\int_{\R^2}|\xi|^{2(s+1)}e^{-(t-\tau) A(\xi)}|\widehat{\Theta.\Theta}|^2d\xi\right)^{\frac{1}{2}}d\tau.
  	\end{align*}
  Moreover, for any $\theta_1,$ $\theta_2\in \mathcal{P}(T)$,we have
  \begin{align*} (1+|\xi|^{2})^{\frac{s}{2}}e^{\frac{t}{2} B(\xi)}|\widehat{\psi(\theta_1)}(t)-\widehat{\psi(\theta_2)}(t)|&\leq C e^T\int_{0}^{t}(1+|\xi|^2)^{\frac{s}{2}}|\xi|e^{-\frac{1}{2}(t-\tau) A(\xi)}e^{\frac{\tau}{2}B(\xi)}|\mathcal{F}({(\theta_1-\theta_2) u_{\theta_1}})|d\tau\\
  	&\hskip0.3cm+C e^{T}\int_{0}^{t}(1+|\xi|^{2})^{\frac{s}{2}}|\xi|e^{-\frac{1}{2}(t-\tau) A(\xi)}e^{\frac{\tau}{2}B(\xi)}|\mathcal{F}({\theta_2 (u_{\theta_1-\theta_2})})|d\tau\\
  	&\leq Ce^T \int_{0}^{t}(1+|\xi|^{2})^{\frac{s}{2}}|\xi|e^{-\frac{1}{2}(t-\tau) A(\xi)}|\mathcal{F}((\Theta_1+\Theta_2).g)|d\tau,
  \end{align*}
  where  $\widehat{\Theta_1}(\xi)=e^{\frac{t}{2} B(\xi)} |\widehat{\theta_1}(\xi)|$, $\widehat{\Theta_2}(\xi)=e^{\frac{t}{2} B(\xi)} |\widehat{\theta_2}(\xi)|$ and $\widehat{g}(\xi)=e^{\frac{t}{2} B(\xi)} |\widehat{\theta_1-\theta_2}(\xi)|$.\\

   Therefore
  	\begin{align*}
  		\|e^{\frac{t}{2} B(D)}(\psi(\theta_1)-\psi(\theta_2))\|_{H^s}&\leq  Ce^{T}\int_{0}^{t}\left(\int_{\R^2}(1+|\xi|^{2})^{s}|\xi|^2e^{-(t-\tau) A(\xi)}|\mathcal{F}((\Theta_1+\Theta_2).g)|^2d\xi\right)^{\frac{1}{2}}d\tau.
\end{align*}
  	The work consists of two parts:\\
  	$\bullet$ If $s\in (2-2\alpha,1)$, inspiring to \eqref{B(u,v)}, we get
  	\begin{align*}
  		&\int_{0}^{t}\left(\int_{\R^2}|\xi|^2e^{-(t-\tau) A(\xi)}|\widehat{\Theta_1.\Theta_2}|^2d\xi\right)^{\frac{1}{2}}d\tau+\int_{0}^{t}\left(\int_{\R^2}|\xi|^{2(s+1)}e^{-(t-\tau) A(\xi)}|\widehat{\Theta_1.\Theta_2}|^2d\xi\right)^{\frac{1}{2}}d\tau\\
  		&\hskip4cm\leq C(\alpha,\beta,s)  \left(T^{\frac{s-2+2\alpha}{2\alpha}}+T^{\frac{s-2+2\beta}{2\beta}}\right)\left\|\Theta_1\right\|_{L^\infty_T(H^s)}\left\|\Theta_2\right\|_{L^\infty_T(H^s)}.
  	\end{align*}
  	The fact that $\theta\in \mathcal{P}(T)$, gives
  	$$\left\|\Theta\right\|_{H^s}=\left\|e^{\frac{t}{2}B(D)}\theta\right\|_{H^s}\leq 2\left\|\theta^0\right\|_{H^s}.$$
  	Therefore
  	\begin{align} \|e^{\frac{t}{2} B(D)}\psi(\theta)\|_{H^s}
  		&\leq e^{T}\|\theta^0\|_{H^s}+4e^{T} C_3(\alpha,\beta,s)  \left(T^{\frac{s-2+2\alpha}{2\alpha}}+T^{\frac{s-2+2\beta}{2\beta}}\right)\left\|\theta^0\right\|_{H^s}^2
  	\end{align}
  and
  	\begin{align*}
	\|e^{\frac{t}{2} B(D)}(\psi(\theta_1)-\psi(\theta_2))\|_{H^s}&\leq e^{T} C_3(\alpha,\beta,s)  \left(T^{\frac{s-2+2\alpha}{2\alpha}}+T^{\frac{s-2+2\beta}{2\beta}}\right)\left\|\Theta_1+\Theta_2\right\|_{L^\infty_T(H^s)}\left\|g\right\|_{L^\infty_T(H^s)}\\
	&\hskip-3cm\leq e^{T}C_3(\alpha,\beta,s)  \left(T^{\frac{s-2+2\alpha}{2\alpha}}+T^{\frac{s-2+2\beta}{2\beta}}\right)\left(\left\|\Theta_1\right\|_{L^\infty_T(H^s)}+\left\|\Theta_2\right\|_{L^\infty_T(H^s)}\right)\left\|e^{\frac{.}{2} B(D)} (\theta_1-\theta_2)\right\|_{L^\infty_T(H^s)}\\
	&\hskip-3cm\leq 4e^{T}C_3(\alpha,\beta,s) \|\theta^0\|_{H^s} \left(T^{\frac{s-2+2\alpha}{2\alpha}}+T^{\frac{s-2+2\beta}{2\beta}}\right)\left\|e^{\frac{.}{2} B(D)} (\theta_1-\theta_2)\right\|_{L^\infty_T(H^s)}.
\end{align*}
  	$\bullet$ If $1\leq s<2$, with the same method for the proof of local existence, we have
  	\begin{align*}
  		&\int_{0}^{t}\left(\int_{\R^2}|\xi|^2e^{-(t-\tau) A(\xi)}|\widehat{\Theta_1. \Theta_2}|^2d\xi\right)^{\frac{1}{2}}d\tau+\int_{0}^{t}\left(\int_{\R^2}|\xi|^{2(s+1)}e^{-(t-\tau) A(\xi)}|\widehat{\Theta_1. \Theta_2}|^2d\xi\right)^{\frac{1}{2}}d\tau\\
  		&\hskip2cm\leq C(\alpha,\beta,s)  \left(T^{\frac{2\alpha-1}{2\alpha}}+T^{\frac{2\beta-1}{2\beta}}+T^{\frac{2\alpha-1}{4\alpha}}+T^{\frac{4\beta-2\alpha-1}{4\beta}}\right)\left\|\Theta_1\right\|_{L^\infty_T(H^s)}\left\|\Theta_2\right\|_{L^\infty_T(H^s)}^2.
  	\end{align*}
  	which imply
  	\begin{align} \|e^{\frac{t}{2} B(D)}\psi(\theta)\|_{H^s}
  		&\leq e^{T}\|\theta^0\|_{H^s}+4e^{T} C_4(\alpha,\beta,s)  \left(T^{\frac{2\alpha-1}{2\alpha}}+T^{\frac{2\beta-1}{2\beta}}+T^{\frac{2\alpha-1}{4\alpha}}+T^{\frac{4\beta-2\alpha-1}{4\beta}}\right)\left\|\theta^0\right\|_{H^s}^2,
  	\end{align}
  and
    	\begin{align*}
  	&	\|	e^{\frac{t}{2}B(D)}(\psi(\theta_1)-\psi(\theta_2))\|_{H^s} \\
  	&	\hskip1cm\leq 4 e^T C_4(\alpha,\beta,s) \|\theta^0\|_{H^s} \left(T^{\frac{2\alpha-1}{2\alpha}}+T^{\frac{2\beta-1}{2\beta}}+T^{\frac{2\alpha-1}{4\alpha}}+T^{\frac{4\beta-2\alpha-1}{4\beta}}\right)\left\|e^{\frac{.}{2}B(D)}(\theta_1-\theta_2)\right\|_{L^\infty_T(H^s)}.
  \end{align*}
  	If we choose $T_1>0$ such that
  	\begin{equation}\label{T}
  		\begin{cases}
  			e^{T_1}<\frac{3}{2},\ 0<T_1\leq T_0\\
  			\\
  			\left(T_1^{\frac{s-2+2\alpha}{2\alpha}}+T_1^{\frac{s-2+2\beta}{2\beta}}\right)e^{T_1}\leq \frac{1}{8C_3(\alpha,\beta,s)\|\theta^0\|_{H^{s}}},\\
  			\\
  			 \left(T_1^{\frac{2\alpha-1}{2\alpha}}+T_1^{\frac{2\beta-1}{2\beta}}+T_1^{\frac{2\alpha-1}{4\alpha}}+T_1^{\frac{4\beta-2\alpha-1}{4\beta}}\right)e^{T_1}\leq \frac{1}{8C_4(\alpha,\beta,s)\|\theta^0\|_{H^{s}}}.\\
  		\end{cases}
  	\end{equation}
  	then, $\psi: \mathcal{B}(T_1)\rightarrow \mathcal{B}(T_1)$  is contracting. The fixed point theorem implies the existence of a unique solution $\theta$ of $(AQG)$ in $C([0,T_1],H^{s}(\R^2))$, moreover  $\psi(\mathcal{P}(T_1)) \subset \mathcal{P}(T_1)$ and it is contracting:
  	Therefore,  The same theorem verify  the existence of a unique solution $\theta'$ of $(AQG)$ in $ \mathcal{P}(T_1)\subset \mathcal{B}(T_1).$\\
  	
  	The uniqueness of solution implies that $\theta=\theta'\in C([0,T_1],H^{s}(\R^2))$ satisfy $$e^{\frac{t}{2}B(D)}\theta\in L^\infty([0,T_1],H^{s}(\R^2)).$$
  	Which completes the proof of Proposition \ref{thm1.3}.
\subsection{Step 3: Regularity $H^2$ of local solution} By the previous steps, there exists
  	$T_0\geq T_1>0$  such that The system $(AQG)$ admits a unique solution
  	$\theta\in C([0,T_0];H^{s}(\R^2))$ satisfy
  	$$
  		\left(t\mapsto e^{t\left(|\partial_{1}|^{\alpha}+|\partial_{2}|^{\beta}\right)}\theta (t)\right)\in L^\infty([0,T_1];H^{s}(\R^2)).
  	$$
  	We want to prove that $\theta\in C((0,T_1),H^2(\R^2))$.\\
  	
  	Let $t_0\in (0,T_1)$, then there exist $\varepsilon>0$ such that $ [t_0-\varepsilon,t_0+\varepsilon]\subset (0,T_1)$. So we have, for  $t\in (t_0-\varepsilon,t_0+\varepsilon)$
  	\begin{align*}
  		\|\theta(t)-\theta(t_0)\|_{H^2}^2&=\int_{\R^2}(1+|\xi|^2)^2 \left|\widehat{\theta}(t,\xi)-\widehat{\theta}(\xi,t_0)\right|^2 d\xi\\
  		&=\int_{\R^2}(1+|\xi|^2)^{\frac{s}{2}} \left|\widehat{\theta}(t,\xi)-\widehat{\theta}(\xi,t_0)\right|(1+|\xi|^2)^{2-\frac{s}{2}} \left|\widehat{\theta}(t,\xi)-\widehat{\theta}(\xi,t_0)\right| d\xi\\
  		&\leq \|\theta(t)-\theta(t_0)\|_{H^s}\|\theta(t)-\theta(t_0)\|_{H^{4-s}}.
  	\end{align*}
  	But, we have
  	\begin{align*}
  		\|\theta(t)-\theta(t_0)\|_{H^{4-s}}^2&= \int_{\R^2}(1+|\xi|^2)^{4-2s} e^{-(t_0-\varepsilon)B(\xi)}(1+|\xi|^2)^{s}e^{(t_0-\varepsilon)B(\xi)} \left|\widehat{\theta}(t,\xi)-\widehat{\theta}(\xi,t_0)\right|^2 d\xi\\
  		&= 2\int_{\R^2}(1+|\xi|^2)^{4-2s} e^{-(t_0-\varepsilon)B(\xi)}(1+|\xi|^2)^{s} \left(e^{tB(\xi)}|\widehat{\theta}(t,\xi)|^2+e^{t_0B(\xi)}|\widehat{\theta}(\xi,t_0)|^2\right) d\xi.
  	\end{align*}
  	Moreover
  	\begin{align*}
  		(1+|\xi|^2)^{4-2s} e^{-(t_0-\varepsilon)B(\xi)}\leq C\left(1+ \frac{1}{(t_0-\varepsilon)^{\frac{8-4s}{\alpha}}}+\frac{1}{(t_0-\varepsilon)^{\frac{8-4s}{\beta}}}\right).
  	\end{align*}
  	Therefore
  	\begin{align*}
  		\|\theta(t)-\theta(t_0)\|_{H^2}^2&
  		\leq C \left(1+ \frac{1}{(t_0-\varepsilon)^{\frac{4-2s}{\alpha}}}+\frac{1}{(t_0-\varepsilon)^{\frac{4-2s}{\beta}}}\right)\|e^{\frac{\centerdot}{2}B(D)}\theta\|_{L^\infty_{T_0}(H^s)}\|\theta(t)-\theta(t_0)\|_{H^{s}}.
  	\end{align*}
  	The fact that $\theta\in C([0,T_0],H^s(\R^2))$, which implies, for any $t_0\in (0,T_0)$
  	\begin{equation}
  		\lim\limits_{t\rightarrow t_0}\|\theta(t)-\theta(t_0)\|_{H^2}=0.
  	\end{equation}
  	Therefore $\theta\in C((0,T_1),H^2(\R^2))$.\\
\subsection{Step 4: Global existence}  	
  	Now, we consider the following system
  		$$(AQG_1)	\left\{\begin{array}{l}
  			\partial_t\gamma+ u_\gamma.\nabla\gamma +|\partial_1|^{2\alpha}\gamma+ |\partial_2|^{2\beta}\gamma=0,\\
  			u_\gamma=\mathcal{R}^\perp \gamma,\\
  			\gamma(0)=\theta\left(\frac{T_0}{2}\right)\in H^2(\R^2).
  		\end{array}\right.$$

By Theorem \ref{thmYE}, there exists a unique global solution of $(AQG_1)$:  $$\gamma \in C(\R^+,H^2(\R^2))\subset C(\R^+,H^s(\R^2)).$$
  	
  	We notice
  	$$
  		\theta^\ast(t)=\left\{\begin{array}{l}
  			\theta(t),\;if\;t\in \left[0,\frac{T_0}{2}\right]\\
  			\gamma(t-\frac{T_0}{2}),\;if\;t\in \left[\frac{T_0}{2},+\infty\right).
  	\end{array}\right.$$
  	By the uniqueness of solution we get $\theta^\ast$ is a global solution of $(AQG)$ in $C(\R^+,H^{s}(\R^2))$, which complete the proof of the main result.
\section{\bf General remarks}
\begin{enumerate}
\item[\textbf{1)}] Adapting the same technique in the third step of the proof, we can prove the Gevrey-class regularity of the solution with $\varepsilon_1\in(0,\alpha)$ instead of $\alpha$ and $\varepsilon_2\in(0,\beta)$ instead of $\beta$, by showing that
$$\left(t\mapsto e^{t(|\partial_{1}|^{\varepsilon_1}+|\partial_{2}|^{\varepsilon_2})}\theta\right)\in L^\infty([0,T_1],H^s(\R^2)).$$
The idea is inspired by the associated linear system:
\begin{equation*}
\left\{\begin{array}{l}
	\partial_t\theta +(|\partial_1|^{2\alpha}+ |\partial_2|^{2\beta})\theta=0,\\
	\theta(x,0)=\theta^0(x).
\end{array}\right.
\end{equation*}
But, the technique used, cannot  be applied to $\varepsilon_1>1$ or $\varepsilon_2>1$, because the inequality
$$|x+y|^a\leq |x|^a+|y|^a;\forall x,y\in \R$$
is false for $a>1$ (see \cite{BJA}).
\item[\textbf{2)}] The fact that (\ref{1.2}) gives the Gevrey-class regularity for the solution is due to the following elementary inequalities:
$$|\xi|^\alpha\leq |\xi_1|^\alpha+|\xi_2|^\alpha\leq |\xi_1|^\alpha+|\xi_2|^\beta+1\leq 2|\xi_2|^\beta+2\leq 2|\xi|^\beta+2 $$
and
$$e^{-T_0}e^{t|\xi|^\alpha}\leq e^{t(|\xi_1|^\alpha+|\xi_2|^\beta)}\leq e^{T_0}e^{2t|\xi|^\beta},\;\forall t\in[0,T_1].$$
Particularly, we get
$$\theta(t)\in H^s_{t,\alpha^{-1}}(\R^2),\;\forall t\in(0,T_1).$$
\item[\textbf{3)}] There are some open problem:
	\begin{itemize}
		\item[$\bullet$] Global existence in Sobolev space $H^{\max\{2-2\alpha,2-2\beta\}}(\R^2)$, when $(\alpha,\beta)\in Y$.
		\item[$\bullet$] Global existence in Sobolev space $H^{s}(\R^2)$, $s\in (\max\{2-2\alpha,2-2\beta\},2),$ when $(\alpha,\beta)\in Y_2\cup Y_3$.
		\item[$\bullet$] The asymptotic study of the global solution of our system in the neighborhood of infinity in the Sobolev space $H^{s}(\R^2)$, $s\in (\max\{2-2\alpha,2-2\beta\},2).$
	\end{itemize}

\end{enumerate}

\end{document}